\documentclass[1p]{elsarticle}
\usepackage{amssymb}
\usepackage{amsfonts}

\newtheorem{theorem}{Theorem}
\newtheorem{conjecture}[theorem]{Conjecture}
\newtheorem{corollary}[theorem]{Corollary}
\newtheorem{lemma}[theorem]{Lemma}

\newproof{pf}{Proof}

\begin{document}
\title{Asymptotically optimal neighbour sum distinguishing total colourings of graphs}

\author{Jakub Przyby{\l}o\fnref{grantJP,MNiSW}}
\ead{jakubprz@agh.edu.pl, phone: 048-12-617-46-38,  fax: 048-12-617-31-65}

\fntext[grantJP]{Supported by the National Science Centre, Poland, grant no. 2014/13/B/ST1/01855.}
\fntext[MNiSW]{Partly supported by the Polish Ministry of Science and Higher Education.}

\address{AGH University of Science and Technology, al. A. Mickiewicza 30, 30-059 Krakow, Poland}

\begin{abstract}
Consider a simple graph $G=(V,E)$ of maximum degree $\Delta$
and its \emph{proper total} colouring $c$
with the elements of the set $\{1,2,\ldots,k\}$.
The colouring $c$ is said to be \emph{neighbour sum distinguishing}
if for every pair of \emph{adjacent} vertices $u$, $v$, we have
$c(u)+\sum_{e\ni u}c(e)\neq c(v)+\sum_{e\ni v}c(e)$.
The least
integer $k$ for which it
exists is
denoted by $\chi''_{\sum}(G)$, hence $\chi''_{\sum}(G) \geq \Delta+1$.
On the other hand, it has been daringly conjectured that
just one more label than presumed in the famous Total Colouring Conjecture
suffices to construct such total colouring $c$,
i.e., that $\chi''_{\sum}(G) \leq \Delta+3$
for all graphs.
We support this inequality by
proving its asymptotic version,
$\chi''_{\sum}(G) \leq (1+o(1))\Delta$.
The major part of the construction
confirming this relays on
a random assignment of colours, where the choice for every edge is biased by so called
\emph{attractors}, randomly assigned to the vertices,
and the probabilistic result of Molloy and Reed on the Total Colouring Conjecture itself.
\end{abstract}


\maketitle

\section{Introduction}

\subsection{Origins}
One of the most elementary facts
we learn in the very first lecture
of a basic combinatorial course
is that every (simple) graph of order at least two
contains a pair of vertices of the same degree.
This datum gave rise to the natural question studied e.g.
by Chartrand, Erd\H{o}s and Oellermann in~\cite{ChartrandErdosOellermann},
on a possible definition of an \emph{irregular graph},
intended as the antonym to the term `regular graph'.
With no convincing individual
solution to the problem, Chartrand et al.~\cite{Chartrand} altered
towards measuring the `irregularity of a graph' instead.
Suppose that given a graph $G=(V,E)$ we wish to construct a multigraph
with pairwise distinct vertex degrees of it by multiplying some of its edges.
The least $k$ so that we are able to achieve such goal
using at most $k$ copies of every edge
is known as the \emph{irregularity strength} of $G$ and denoted by $s(G)$, see~\cite{Chartrand}.
Alternatively, one may consider (not necessarily proper) edge colourings
$c:E\to\{1,2,\ldots,k\}$ with $\sum_{e\ni u}c(e)\neq\sum_{e\ni v}c(e)$
for every pair of distinct vertices $u,v\in V$. Then the least
$k$ which permits defining a colouring $c$ with this feature
equals $s(G)$.
It is straightforward to notice that $s(G)$ is well defined for all graphs
containing no isolated edges and at most one isolated vertex.
The irregularity strength was studied in numerous papers, e.g.~\cite{Aigner,Bohman_Kravitz,Lazebnik,Dinitz,Faudree,Frieze,KalKarPf,Lehel,MajerskiPrzybylo2,Nierhoff,Przybylo,irreg_str2},
and was
the cornerstone of the later \emph{additive graph labelings},
or more generally -- \emph{vertex distinguishing graph colourings}.
Many consequential and related graph parameters have been studied ever since its development.
These are associated with
quite a few new interesting proving methods
and deep results reaching beyond this particular field.

One of the
most closely related successor
of the irregularity strength
 was the problem
the only novelty concerning which asserted distinguishing only adjacent vertices
with their corresponding sums, i.e., its concern was the least $k$ so that
a colouring $c:E\to\{1,2,\ldots,k\}$ existed with $\sum_{e\ni u}c(e)\neq\sum_{e\ni v}c(e)$
for \emph{every edge} $uv\in E$.
In fact Karo\'nski \L uczak and Thomason~\cite{123KLT} conjectured that $k=3$ is sufficient
for every connected graph of order at least $3$. Their presumption
is now well known and studied in the literature under a common name
of the \emph{1--2--3 Conjecture}, see e.g.~\cite{Louigi30,Louigi}.
Thus far it is known that $k=5$ suffice, \cite{KalKarPf_123}.
Another noted conjecture of the field comes from~\cite{Zhang} and is commonly
referred to as the \emph{Zhang's Conjecture}.
This time the edge colouring $c$ is required to be \emph{proper}, while the adjacent vertices
are supposed to differ only in the \emph{sets} of their incident colours,
    rather than sums.
Zhang, Liu and Wang~\cite{Zhang} conjectured
that only slightly more than required in a proper edge colouring, i.e., $k=\Delta+2$ colours
are always sufficient
for designing such $c$ for every connected graph of order at least three,
unless it is the cycle $C_5$.
The corresponding parameter, known as the \emph{neighbour set distinguishing index}, see e.g.~\cite{FlandrinMPSW},
or the \emph{adjacent strong chromatic index}, see~\cite{Zhang},
and denoted by $\chi'_a(G)$,
came by an almost optimal upper bound delivered by Hatami~\cite{Hatami}.
\begin{theorem}
\label{HatamiTh}
If $G$ is a graph with no isolated edges and with maximum degree $\Delta>10^{20}$, then $\chi'_a(G)\leq \Delta+300$.
\end{theorem}
The proof of this was based on a multistage probabilistic argument.
See also~\cite{Akbari,BalGLS,HocqMont,Hornak_planar}
for some results and other notations concerning this parameter.
Proper colourings were the setting for yet another graph invariant, related directly with the
1--2--3 Conjecture.
The least $k$ so that there is a \emph{proper} edge colouring $c:E\to\{1,2,\ldots,k\}$
such that $\sum_{e\ni u}c(e)\neq\sum_{e\ni v}c(e)$
for every edge $uv\in E$ is called the \emph{neighbour sum distinguishing index}
and denoted by $\chi'_{\sum}(G)$, see~\cite{FlandrinMPSW},
and~\cite{DongWang_planar,Przybylo_asym_optim,Przybylo_CN_1,Przybylo_CN_2}
for some of the results concerning this parameter.
Note also that $\chi'_a(G)\leq \chi'_{\sum}(G)$, and
though the probabilistic method is much more unwieldy in sum setting
(due to the concentration of a sum of independent random variables with uniform distribution, and many other reasons) its application yielded the following upper bound, asymptotically
equivalent to the one proved by Hatami in the case of sets in Theorem~\ref{HatamiTh}.
\begin{theorem}[\cite{Przybylo_asym_optim}]
If $G$ is a connected graph of maximum degree $\Delta\geq 2$, then $\chi'_{\sum}(G)\leq (1+o(1))\Delta$.
\end{theorem}
These problems make up a fair share
of the list of the central issues of the discipline,
followed by multiple their corresponding variants, including e.g. list versions of these.
See~\cite{Aldred,BarGrNiw,Daltonisci_1,PrzybyloWozniakChoos,WongZhu23Choos,WongZhuChoos}
for a few interesting and
influential examples.

\subsection{Total Colouring Conjecture}
Probabilistic approach was also used to obtain a big breakthrough in a
famous colouring problem which has eluded mathematicians for $40$ years.
A \emph{proper} \emph{total colouring} of $G$
is a colouring of its vertices and edges so that no two adjacent vertices get the same colour,
no two incident edges get the same colour, and no edge gets the same colour as one of its endpoints.
The least number of colours in such a colouring is called the \emph{total chromatic number} of $G$ and is denoted by $\chi''(G)$, hence $\chi''(G)\geq \Delta(G)+1$.
The following conjecture was posed by Vizing~\cite{Vizing2} and independently by Behzad~\cite{Behzad}.
\begin{conjecture}[The Total Colouring Conjecture]
For every graph $G$, $\chi''(G)$ $\leq$ $\Delta(G)$ $+2$.
\end{conjecture}
An asymptotic confirmation of this was delivered by Molloy and Reed \cite{MolloyReedTotal}, who designed a complex probabilistic argument implying that $\chi''(G)\leq\Delta(G)+const.$
\begin{theorem}[\cite{MolloyReedTotal}]\label{MolloyReedTh}
There exists a constant $\Delta_0$ such that for every graph $G$ with $\Delta(G)\geq\Delta_0$, $\chi''(G)\leq\Delta(G)+10^{26}$.
\end{theorem}
They also remarked that being more careful and adding
a few (obscuring the clarity of presentation) additional intricacies
they could prove the same result with the constant $10^{26}$
exchanged with $500$ or even slightly lower number.
Since obviously $\chi''(G)\leq 2\Delta(G)+1$, we thus have:
\begin{corollary}\label{MolloReedCor}
There is a constant $C$
such that $\chi''(G)\leq \Delta(G)+C$ for every graph $G$.
\end{corollary}

\subsection{Main Objective and Tools}

Given any proper total colouring $c:V\cup E\to\{1,2,\ldots,k\}$ of $G$ and $v\in V$,
the sum
$$
s_c(v):=c(v)+\sum_{u\in N(v)}c(uv)
$$
shall be called the \emph{weighted degree} of the vertex $v$.
When it causes no confusion, we shall also write $s(v)$ instead of $s_c(v)$.
The least $k$
which permits constructing such $c$ that attributes distinct
weighted degrees to the adjacent vertices in $G$
is called its
\emph{total neighbour sum distinguishing number} and denoted $\chi''_{\sum}(G)$, see~\cite{PilsniakWozniak_total}.
Other results concerning this graph invariant can be found
in~\cite{total_sum_planar,LiLiuWang_sum_total_K_4,PilsniakWozniak_total},
and in~\cite{Przybylo_CN_3}, where it has been proved that $\chi''_{\sum}(G)\leq \Delta(G)+ \lceil\frac{5}{3}{\rm col}(G)\rceil$.
Though constructing
a proper total colouring using $\Delta$ plus only a few additional
colours is difficult itself, as the long history of the Total Colouring Conjecture exemplifies,
we shall prove that \emph{asymptotically} this many are sufficient to
find one satisfying even our additional condition that $s(u)\neq s(v)$ for every edge $uv\in E$.
In particular we shall prove that $\chi''_{\sum}(G)\leq (1+o(1))\Delta$ for all graphs,
see Theorem~\ref{main_JP_Th_as} below. This asymptotically confirms the following conjecture of Pil\'sniak and Wo\'zniak.
\begin{conjecture}[\cite{PilsniakWozniak_total}]\label{PilsniakWozniakCon}
For every graph $G$, $\chi''_{\sum}(G)\leq \Delta(G)+3$.
\end{conjecture}

It is also worth mentioning here that with the requirement of the \emph{properness} of
the total colourings investigated \emph{skipped}, the correspondent of this problem is known as
the \emph{1--2 Conjecture}, as it is presumed within it that just labels 1 and 2
are always sufficient for designing appropriate (not necessarily proper) total colourings
(so that $s(u)\neq s(v)$ for $uv\in E$),
see e.g.~\cite{Kalkowski12,12Conjecture}.
One may also confront
our key result -- Theorem~\ref{main_JP_Th_as} -- concerning the
total neighbour sum distinguishing number
with those devoted to its
less restrictive counterpart focused on sets
rather than sums, see e.g.~\cite{CokerJohanson,Zhang_total}.

Our approach, similarly as in~\cite{Przybylo_asym_optim}, shall consist of two main parts.
First we shall generate some total colouring using a precisely designed random process,
in which the drawing rules for every edge shall be biased by
so called attractors.
After this careful probabilistic construction the labels of most of the edges and vertices
shall be close to their final values, while weighted degrees shall already be
quite well scattered.
In the second part, the colours of (almost) all edges and vertices
shall be slightly modified so that the total colouring obtained is proper and neighbour sum distinguishing.
More detailed idea of our proof is contained in the next section.
Our approach is however completely different from the one used by Hatami~\cite{Hatami}
and Coker and Johanson~\cite{CokerJohanson} to deal with the set case,
as it was impossible to implant even the main threads of their
ideas concerning sets
into the much less hospitable ground of sums.
To get a rough feeling
 why the sum setting is significantly more troublesome for application of the probabilistic method,
one may e.g. consider a (random) total colouring of a $\Delta$-regular graph
with the colours $1,2,\ldots,\Delta+K$.
Then the number of different $(\Delta+1)$-element \emph{subsets} of the set of colours
$\{1,2,\ldots,\Delta+K\}$, where $\Delta$ is large and $K$ is some constant
(e.g., $K=300$, as in~\cite{Hatami}),
is of order $\Delta^{K}$, while the same subsets generate no more than roughly $K\Delta$
distinct \emph{sums}, whose distribution is -- unlike for sets -- concentrated around
their expected value,
not to mention that in case of
sets we need only to distinguish the neighbours of the same degrees,
contrary to the sum case.

Throughout the paper we shall use several times two classical
tools of the probabilistic method,
the Lov\'asz Local Lemma, see e.g.~\cite{AlonSpencer}, and the Chernoff Bound, see e.g.~\cite{MolloyReed}.
\begin{theorem}[\textbf{The Local Lemma; Symmetric case}]\label{LLL-symmetric}
Let $A_1,A_2,\ldots,A_n$ be events in an arbitrary probability space.
Suppose that each event $A_i$ is mutually independent of a set of all the other
events $A_j$ but at most $D$, and that ${\rm \emph{\textbf{Pr}}}(A_i)\leq p$ for all $1\leq i \leq n$. If
$$ ep(D+1) \leq 1,$$
then $ {\rm \emph{\textbf{Pr}}}\left(\bigcap_{i=1}^n\overline{A_i}\right)>0$.
\end{theorem}
\begin{theorem}[\textbf{Chernoff Bound}]\label{ChernofBoundTh}
For any $0\leq t\leq np$:
$${\rm\emph{\textbf{Pr}}}(|{\rm BIN}(n,p)-np|>t)<2e^{-\frac{t^2}{3np}},$$
where ${\rm BIN}(n,p)$ is the sum of $n$ independent variables, each equal to $1$ with probability $p$ and $0$ otherwise.
\end{theorem}

\section{General Proof Idea}
Assume that a graph $G=(V,E)$ has a large maximum degree $\Delta$, i.e., large enough
so that a few explicit inequalities involved in the proof hold.
We shall in fact prove that $\chi''_{\sum}(G) \leq \Delta+O(\Delta^{\frac{5}{6}}\ln^\frac{1}{6}\Delta)$.
For this aim we shall carry out a carefully designed random experiment
providing, with positive probability, a total colouring with several useful features,
e.g., that the colours are nicely scattered, i.e., the vertices
adjacent to every vertex (of large degree) consist of almost equally-sized
subsets of elements in each of the colours used,
and the same holds for the edges incident with each vertex (of large degree).
The edge colours shall not be chosen uniformly though.
In fact we first randomly pick an auxiliary colour
$c_2(e)\in\{1,2,\ldots,\lceil\Delta^{\frac{1}{3}}/\ln^\frac{1}{3}\Delta\rceil\}$
for every edge $e\in E$,
and an auxiliary colour
$c_1(v)\in\{1,2,\ldots,\lceil\Delta^\frac{1}{6}/\ln^\frac{1}{6}\Delta\rceil\}$
-- chosen randomly from a smaller range of possibilities -- for every $v\in V$,
and only just then define the first approximation of our desired colouring by setting
$c_3(uv)=c_1(u)+c_2(uv)+c_1(v)$ for every edge $uv\in E$
(while the values of $c_3(v)$ are chosen randomly and uniformly for the vertices $v\in V$).
Note that such construction of the total colouring $c_3$,
and especially the definition of edge colours within it,
has a consequential influence on the expected values of the weighted degrees of all vertices,
since for every vertex $v$ of (large) degree $d$,
the value $c_1(v)$ shall be counted $d$ times in its weighted degree $s(v)$.
Thus
the quantity $c_1(v)$ plays the role
of an \emph{attractor}.
It pulls
the weighted degree of $v$ towards a quantity, say $S(v)$,
dependent exclusively on $c_1(v)$, $d(v)$ (and $\Delta$).
The relatively even distributions of the colours $c_1(u)$ and $c_2(uv)$
around every vertex
shall guarantee that the weighted degree of $v$ is `close' to that $S(v)$.
Intuitively, since the neighbours $u$ of $v$
shall have well distributed attractors, $c_1(u)$,
only few of them should
have $S(u)$, hence also the weighted degree, close to $S(v)$,
at least among the neighbours of $v$ of similar degree as $v$.
(Distinguishing weighted degrees of vertices with small degrees
from those of their neighbours shall occur straightforward at the end of our construction.)
We shall in fact prove that our colourings might be chosen so that
for a partition of the positive reals
into relatively short intervals $I_{0},I_{1},I_{2},\ldots$,
for every $j\in\{0,1,2,\ldots\}$ only for a limited number of neighbours $u$ of
any given vertex $v$ of large degree, we have $S(u)\in I_{j}$.

Note that by the limited range of labels allowed within the colourings $c_1$ and $c_2$,
our total colouring $c_3$ cannot be proper. We shall thus recalibrate
it by (temporarily) setting $c_t(e)=B\cdot c_3(e)$ and $c_t(v)=B\cdot c_3(v)$,
where the constant $B$ shall be chosen so that $c_t(e),c_t(v)\leq \Delta+O(\Delta^{\frac{5}{6}}\ln^\frac{1}{6}\Delta)$
for every $e\in E$ and $v\in V$.
Such spreading of the colours of edges and vertices
shall enable then introducing small changes to these
in order to eliminate adjacency of the same colours.
This shall be achieved
via (multiple) application of Theorem~\ref{MolloyReedTh} of Molloy and Reed
to redistribute
the colours in each colour class,
i.e., in each subgraph induced by the edges $e$ with the same colour $c_3(e)$
assigned, separately
(where no colour conflict shall be possible between
objects from different colour classes within our construction).
As the vertices shall retain their resulting colours attributed
to them while analyzing the subgraph corresponding to their colours $c_3$,
to guarantee distinction of the final colours for adjacent vertices,
for every edge $uv\in E$ with $c_3(u)=c_3(v)$ we shall additionally
need that $c_3(uv)=c_3(u)$.
This requires some small technical modifications in the
priory outlined construction of $c_3$, causing
minor deviations from its definition for few of the edges.
All these shifts and minor exceptions
shall not have a significant influence
on the weighted degrees of the vertices of large degree, though.
Thus at the end of the construction we shall be able to use `a few' extra integers
to recolour edges of some random sparse subgraph of $G$
so that neighbours of large degree have pairwise distinct weighted degrees.
By modifying vertex colours we shall then distinguish all vertices
of small degrees from their neighbours.

All other technical details
are included in section~\ref{section_with_main_result}, while the main
probabilistic lemma, proving the existence of auxiliary colourings of the vertices and edges with the mentioned features follows
in the next section.

\section{Main Probabilistic Lemma}
\begin{lemma}\label{main_probabilistic_lemma}
Let $G=(V,E)$ be a graph of maximum degree
$\Delta$.
For any positive integer $\alpha$, let
\begin{eqnarray}
I_{\alpha}&:=&\left((\alpha-1) \frac{\Delta^{\frac{5}{3}}\ln^{\frac{1}{3}}\Delta}{3}, \alpha\frac{\Delta^{\frac{5}{3}}\ln^{\frac{1}{3}}\Delta}{3}\right].\label{I_d_alpha}
\end{eqnarray}
If $\Delta$ is sufficiently large, then there exist colourings
$$c_1:V\to\left\{1,2,\ldots,\left\lceil\frac{\Delta^{\frac{1}{6}}}{\ln^\frac{1}{6}\Delta}\right\rceil\right\},$$
$$c_2:E\to\left\{1,2,\ldots,\left\lceil\frac{\Delta^{\frac{1}{3}}}{\ln^\frac{1}{3}\Delta}\right\rceil\right\}$$
and
$$c_3:V\cup E\to\left\{1,2,\ldots,2\left\lceil\frac{\Delta^{\frac{1}{6}}}{\ln^\frac{1}{6}\Delta}\right\rceil +\left\lceil\frac{\Delta^{\frac{1}{3}}}{\ln^\frac{1}{3}\Delta}\right\rceil\right\}$$
such that if
for every vertex $v\in V$ of degree $d$, $0\leq d\leq \Delta$,
\begin{equation}\label{S_v_in_lemma}
S(v)=\left(\left\lceil\Delta^{\frac{2}{3}}\ln^{\frac{1}{3}}\Delta\right\rceil
+6\left\lceil\Delta^{\frac{1}{3}}\ln^{\frac{2}{3}}\Delta\right\rceil\right)dc_1(v) +R(d,\Delta),
\end{equation}
where $R(d,\Delta)$ is any given function of $d$ and $\Delta$ for which $S(v)\leq\Delta^2$, then for every vertex $v$ of degree $d$:
\begin{itemize}
  \item[(I)] if $d\geq \frac{\Delta}{3}$,
  then for every $c_1^*\in\{1,\ldots,\lceil\Delta^{\frac{1}{6}}/\ln^\frac{1}{6}\Delta\rceil\}$,
      the number of neighbours $u$ of $v$ with $c_1(u)=c_1^*$
      equals $\frac{d}{\lceil\Delta^{\frac{1}{6}}\ln^{-\frac{1}{6}}\Delta\rceil}+f_{1,c_1^*}(v)$,
      where
      $$|f_{1,c_1^*}(v)|\leq \Delta^\frac{1}{2}
      ;$$
  \item[(II)] if $d\geq \frac{\Delta}{3}$,
  then for every $c_2^*\in\{1,\ldots,\lceil\Delta^{\frac{1}{3}}/\ln^\frac{1}{3}\Delta\rceil\}$,
      the number of neighbours $u$ of $v$ with $c_2(uv)=c_2^*$
      equals $\frac{d}{\lceil\Delta^{\frac{1}{3}}\ln^{-\frac{1}{3}}\Delta\rceil}+f_{2,c_2^*}(v)$,
      where
      $$|f_{2,c_2^*}(v)|\leq 3\Delta^{\frac{1}{3}}\ln^{\frac{2}{3}}\Delta;$$
  \item[(III)]
      $c_3(uv)=c_1(u)+c_1(v)+c_2(uv)$ for
      at least $d-(3+o(1))\Delta^\frac{2}{3}\ln^\frac{1}{3}\Delta$ neighbours $u$ of $v$;
  \item[(IV)] for every $c_3^*\in\{1,2,\ldots,2\lceil\Delta^{\frac{1}{6}}/\ln^\frac{1}{6}\Delta\rceil+\lceil\Delta^{\frac{1}{3}}/\ln^\frac{1}{3}\Delta\rceil\}$,
      the number of edges
      $uv$ (incident with $v$) for which $
      c_3(uv)=c_3^*$
      does not exceed
      $\Delta^\frac{2}{3}\ln^\frac{1}{3}\Delta + 5\Delta^\frac{1}{3}\ln^\frac{2}{3}\Delta$;
  \item[(V)]
  for every neighbour $u$ of $v$, if $c_3(u)=c_3(v)$, then also $c_3(uv)=c_3(v)$ ($=c_3(u)$);
  \item[(VI)] if
      $d\geq\frac{\Delta}{3}$,
      then for every integer $\alpha > 0$, the number of neighbours $u$ of $v$ with
      $d(u)\geq \frac{\Delta}{3}$ and
      $S(u)\in I_{\alpha}$ does not exceed
      $$\Delta^{\frac{5}{6}}\ln^{\frac{1}{6}}\Delta+ \Delta^{\frac{1}{2}}.$$
\end{itemize}
\end{lemma}

\begin{pf}
We construct our colourings in two steps.
First $c_1(v),c_2(e)$ and $c_3(v)$ shall be chosen randomly for $v\in V$, $e\in E$, so that
\emph{(I),(II),(VI)} from the thesis hold, and so that
\begin{itemize}
\item[$(1^\circ)$]
for every $v\in V$ and every
$c_3^*\in\{1,2,\ldots,2\lceil\Delta^{\frac{1}{6}}/\ln^\frac{1}{6}\Delta\rceil+\lceil\Delta^{\frac{1}{3}}/\ln^\frac{1}{3}\Delta\rceil\}$,
the number of edges
$uv$ (incident with $v$) for which $c_1(u)+c_1(v)+c_2(uv)=c_3^*$
does not exceed
$\Delta^\frac{2}{3}\ln^\frac{1}{3}\Delta + 3\Delta^\frac{1}{3}\ln^\frac{2}{3}\Delta$, while
the number of edges
$uv$ for which  $c_1(u)+c_1(v)+c_2(uv)=c_3(u)$ or $c_1(u)+c_1(v)+c_2(uv)=c_3(v)$
does not exceed
$2\Delta^\frac{2}{3}\ln^\frac{1}{3}\Delta + 5\Delta^\frac{1}{3}\ln^\frac{2}{3}\Delta$, and
\item[$(2^\circ)$] for every $v\in V$ and every
$c_3^*\in\{1,2,\ldots,2\lceil\Delta^{\frac{1}{6}}/\ln^\frac{1}{6}\Delta\rceil+\lceil\Delta^{\frac{1}{3}}/\ln^\frac{1}{3}\Delta\rceil\}$,
the number of neighbours
$u$ of $v$ with $c_3(u)=c_3^*$
does not exceed
$\Delta^\frac{2}{3}\ln^\frac{1}{3}\Delta
 + 3\Delta^\frac{1}{3}\ln^\frac{2}{3}\Delta$.
\end{itemize}
(Though in fact the formulation of $(2^\circ)$ is much stronger than required for the sake of its further
application, its proof has a simpler notation in this form,
which at the same time does not alter our final result.)
To achieve this we carry out three independent random experiments.
First for every vertex $v\in V$ we randomly and independently choose an integer  $c_1(v)\in\{1,2,\ldots,\lceil\Delta^\frac{1}{6}/\ln^\frac{1}{6}\Delta\rceil\}$ (each with equal probability).
Second, for every edge $e\in E$ we randomly and independently choose an integer $c_2(e)\in\{1,2,\ldots,\lceil\Delta^{\frac{1}{3}}/\ln^\frac{1}{3}\Delta\rceil\}$,
and third, for every vertex $v\in V$ we randomly and independently choose one more integer  $c_3(v)\in\{1,2,\ldots,\lceil\Delta^\frac{1}{3}/\ln^\frac{1}{3}\Delta\rceil\}$.
In the following, whenever needed, we shall assume that $\Delta$ is sufficiently large.

Let us first bound the probability that \emph{(I),(II)} or $(2^\circ)$ fails to
hold for some given vertex.
For every $v\in V$ of degree $d$, let $A_v$ denote the event that $d\geq \Delta/3$ and for at least one integer
$c_1^*\in\{1,\ldots,\lceil\Delta^{\frac{1}{6}}/\ln^\frac{1}{6}\Delta\rceil\}$,
the number of neighbours $u$ of $v$ with $c_1(u)=c_1^*$ is outside the
range postulated in \emph{(I)} above,
let $B_v$ denote the event that $d\geq \Delta/3$
and for at least one integer
$c_2^*\in\{1,\ldots,\lceil\Delta^{\frac{1}{3}}/\ln^\frac{1}{3}\Delta\rceil\}$,
the number of neighbours $u$ of $v$ with $c_2(uv)=c_2^*$ is outside the
range postulated in \emph{(II)},
and analogously,
let $C_v$ denote the event that
for at least one integer
$c_3^*\in\{1,\ldots,\lceil\Delta^{\frac{1}{3}}/\ln^\frac{1}{3}\Delta\rceil\}$,
the number of neighbours $u$ of $v$ with $c_3(u)=c_3^*$ is outside the
range postulated in $(2^\circ)$.
Let $X_{v,c_1^*}$, $Y_{v,c_2^*}$ and $Z_{v,c_3^*}$ be the random variables of the numbers of neighbours $u$ of $v$ with $c_1(u)=c_1^*$, $c_2(uv)=c_2^*$ and $c_3(u)=c_3^*$, respectively, $c_1^*\in\{1,\ldots,\lceil\Delta^{\frac{1}{6}}/\ln^\frac{1}{6}\Delta\rceil\}$, $c_2^*,c_3^*\in\{1,\ldots,\lceil\Delta^{\frac{1}{3}}/\ln^\frac{1}{3}\Delta\rceil\}$.
Then $X_{v,c_1^*}\sim {\rm BIN}(d,1/\lceil\Delta^{\frac{1}{6}}\ln^{-\frac{1}{6}}\Delta\rceil)$ and $Y_{v,c_2^*},Z_{v,c_3^*}\sim {\rm BIN}(d,1/\lceil\Delta^{\frac{1}{3}}\ln^{-\frac{1}{3}}\Delta\rceil)$,
and hence by the Chernoff Bound, if $d\geq \Delta/3$, then
\begin{eqnarray*}{\rm \textbf{Pr}}\left(\left|X_{v,c_1^*}-\frac{d}{\lceil\Delta^{\frac{1}{6}}\ln^{-\frac{1}{6}}\Delta\rceil}\right| >\Delta^\frac{1}{2}\right)
&<&2e^{-\frac{\Delta\lceil\Delta^{\frac{1}{6}}\ln^{-\frac{1}{6}}\Delta\rceil}{3d}}\\
&\leq&
2e^{-3\ln\Delta} = \frac{2}{\Delta^3}
\end{eqnarray*}
and analogously,
\begin{eqnarray*}
{\rm \textbf{Pr}}\left(\left|Y_{v,c_2^*}-\frac{d}{\lceil\Delta^{\frac{1}{3}}\ln^{-\frac{1}{3}}\Delta\rceil}\right| >3\Delta^{\frac{1}{3}}\ln^{\frac{2}{3}}\Delta\right)
&<& 2e^{-3\Delta^{\frac{2}{3}}\ln^{\frac{4}{3}}\Delta\frac{\lceil\Delta^{\frac{1}{3}}\ln^{-\frac{1}{3}}\Delta\rceil}{d}}\\
&\leq& 2e^{-3\frac{\Delta}{d}\ln\Delta} \leq 2e^{-3\ln\Delta} = \frac{2}{\Delta^3},
\end{eqnarray*}
and for every $d$,
\begin{eqnarray}
&&{\rm \textbf{Pr}}\left(Z_{v,c_3^*}> \Delta^\frac{2}{3}\ln^\frac{1}{3}\Delta + 3\Delta^\frac{1}{3}\ln^\frac{2}{3}\Delta\right)\nonumber\\
&\leq& {\rm \textbf{Pr}}\left(BIN\left(\Delta,\frac{1}{\lceil\Delta^{\frac{1}{3}}\ln^{-\frac{1}{3}}\Delta\rceil}\right)
> \Delta^\frac{2}{3}\ln^\frac{1}{3}\Delta + 3\Delta^\frac{1}{3}\ln^\frac{2}{3}\Delta\right)\nonumber\\
&\leq&
{\rm\textbf{Pr}}\left(\left|BIN\left(\Delta,\frac{1}{\lceil\Delta^{\frac{1}{3}}\ln^{-\frac{1}{3}}\Delta\rceil}\right)
-\frac{\Delta}{\lceil\Delta^{\frac{1}{3}}\ln^{-\frac{1}{3}}\Delta\rceil}\right| >3\Delta^{\frac{1}{3}}\ln^{\frac{2}{3}}\Delta\right)\nonumber\\
&<& 2e^{-3\Delta^{\frac{2}{3}}\ln^{\frac{4}{3}}\Delta\frac{\lceil\Delta^{\frac{1}{3}}\ln^{-\frac{1}{3}}\Delta\rceil}{\Delta}}
\leq
2e^{-3\ln\Delta} = \frac{2}{\Delta^3},\label{Z_vc_last_ineq}
\end{eqnarray}
Consequently,
\begin{equation}\label{bound_for_A_v}
{\rm \textbf{Pr}}(A_v) \leq \left\lceil\frac{\Delta^{\frac{1}{6}}}{\ln^\frac{1}{6}\Delta}\right\rceil \cdot  \frac{2}{\Delta^3} < \frac{1}{\Delta^2\cdot\Delta^\frac{5}{6}},
\end{equation}
\begin{equation}\label{bound_for_B_v}
{\rm \textbf{Pr}}(B_v) \leq\left\lceil\frac{\Delta^{\frac{1}{3}}}{\ln^\frac{1}{3}\Delta}\right\rceil \cdot  \frac{2}{\Delta^3} < \frac{1}{\Delta^2\cdot\Delta^\frac{2}{3}},
\end{equation}
and
\begin{equation}\label{bound_for_C_v_NEW}
{\rm \textbf{Pr}}(C_v) \leq\left\lceil\frac{\Delta^{\frac{1}{3}}}{\ln^\frac{1}{3}\Delta}\right\rceil
 \cdot  \frac{2}{\Delta^3} < \frac{1}{\Delta^2\cdot\Delta^\frac{2}{3}}.
\end{equation}

Now let us consider our requirement $(1^\circ)$.
For a vertex $v\in V$ of degree $d$, let $D_v$ denote the event that there exists an integer $c_3^*$
such that the number of edges $uv$ (incident with $v$) for which $c_1(u)+c_1(v)+c_2(uv)=c_3^*$
exceeds
$\Delta^\frac{2}{3}\ln^\frac{1}{3}\Delta + 3\Delta^\frac{1}{3}\ln^\frac{2}{3}\Delta$,
and let $E_v$ denote the event that the number of
edges
$uv$ for which $c_1(u)+c_1(v)+c_2(uv)=c_3(u)$ or $c_1(u)+c_1(v)+c_2(uv)=c_3(v)$
exceeds
$2\Delta^\frac{2}{3}\ln^\frac{1}{3}\Delta + 5\Delta^\frac{1}{3}\ln^\frac{2}{3}\Delta$.
For every integer
$c_3^*\in\{1,2,\ldots,2\lceil\Delta^{\frac{1}{6}}/\ln^\frac{1}{6}\Delta\rceil+\lceil\Delta^{\frac{1}{3}}/\ln^\frac{1}{3}\Delta\rceil\}$,
let $T_{v,c_3^*}$ be the random variable of the number of edges $uv$ (incident with $v$) with $c_1(u)+c_1(v)+c_2(uv)=c_3^*$,
and let $Q_{v}$ be the random variable of the number of edges $uv$ with
$c_1(u)+c_1(v)+c_2(uv)=c_3(u)$ or $c_1(u)+c_1(v)+c_2(uv)=c_3(v)$.
Then given $v\in V$ and $c_3^*\in\{1,\ldots,2\lceil\Delta^{\frac{1}{6}}/\ln^\frac{1}{6}\Delta\rceil+\lceil\Delta^{\frac{1}{3}}/\ln^\frac{1}{3}\Delta\rceil\}$, for any fixed vertex colourings $c_1$, $c_3$ and $u\in N(v)$, the probability that $c_1(u)+c_1(v)+c_2(uv)=c_3^*$ (i.e., that $c_2(uv)=c_3^*-c_1(u)-c_1(v)$)
equals at most $1/\lceil\Delta^{\frac{1}{3}}\ln^{-\frac{1}{3}}\Delta\rceil$,
while the probability that $c_1(u)+c_1(v)+c_2(uv)=c_3(u)$ or $c_1(u)+c_1(v)+c_2(uv)=c_3(v)$
(i.e., that $c_2(uv)=c_3(u)-c_1(u)-c_1(v)$ or $c_2(uv)=c_3(v)-c_1(u)-c_1(v)$)
equals at most $2/\lceil\Delta^{\frac{1}{3}}\ln^{-\frac{1}{3}}\Delta\rceil$.
Therefore, since the choices of $c_2(uv)$ for all neighbours $u$ of $v$ are independent,
by the Chernoff Bound we obtain (to be strict, we should have first written below
the conditional probabilities with respect to some fixed colourings $c_1$, $c_3$,
but since the choices for $c_1$, $c_3$ and $c_2$ are independent and we would have obtained the same upper bound regardless of the choices for $c_1$ and $c_3$,
then the application of the total probability would yield what follows):
\begin{eqnarray*}
&&{\rm \textbf{Pr}}\left(T_{v,c_3^*}> \Delta^\frac{2}{3}\ln^\frac{1}{3}\Delta + 3\Delta^\frac{1}{3}\ln^\frac{2}{3}\Delta\right)\\
&\leq& {\rm \textbf{Pr}}\left(BIN\left(\Delta,\frac{1}{\lceil\Delta^{\frac{1}{3}}\ln^{-\frac{1}{3}}\Delta\rceil}\right)
> \Delta^\frac{2}{3}\ln^\frac{1}{3}\Delta + 3\Delta^\frac{1}{3}\ln^\frac{2}{3}\Delta\right)
< \frac{2}{\Delta^3}
\end{eqnarray*}
(see~(\ref{Z_vc_last_ineq})), and:
\begin{eqnarray*}
&&{\rm \textbf{Pr}}\left(Q_{v}> 2\Delta^\frac{2}{3}\ln^\frac{1}{3}\Delta + 5\Delta^\frac{1}{3}\ln^\frac{2}{3}\Delta\right)\\
&\leq& {\rm \textbf{Pr}}\left(BIN\left(\Delta,\frac{2}{\lceil\Delta^{\frac{1}{3}}\ln^{-\frac{1}{3}}\Delta\rceil}\right)
> 2\Delta^\frac{2}{3}\ln^\frac{1}{3}\Delta + 5\Delta^\frac{1}{3}\ln^\frac{2}{3}\Delta\right)\\
&\leq&
{\rm\textbf{Pr}}\left(\left|BIN\left(\Delta,\frac{2}{\lceil\Delta^{\frac{1}{3}}\ln^{-\frac{1}{3}}\Delta\rceil}\right)
-\frac{2\Delta}{\lceil\Delta^{\frac{1}{3}}\ln^{-\frac{1}{3}}\Delta\rceil}\right| >5\Delta^{\frac{1}{3}}\ln^{\frac{2}{3}}\Delta\right)\\
&<& 2e^{-25\Delta^{\frac{2}{3}}\ln^{\frac{4}{3}}\Delta\frac{1}{3}
\frac{\lceil\Delta^{\frac{1}{3}}\ln^{-\frac{1}{3}}\Delta\rceil}{2\Delta}}
\leq
2e^{-3\ln\Delta} = \frac{2}{\Delta^3},
\end{eqnarray*}
Therefore,
\begin{equation}\label{bound_for_D_v_NEW}
{\rm \textbf{Pr}}(D_v)<\left(2\left\lceil\frac{\Delta^{\frac{1}{6}}}{\ln^\frac{1}{6}\Delta}\right\rceil
+\left\lceil\frac{\Delta^{\frac{1}{3}}}{\ln^\frac{1}{3}\Delta}\right\rceil\right)\frac{2}{\Delta^3}\leq \frac{1}{\Delta^2\cdot\Delta^\frac{2}{3}}
\end{equation}
and
\begin{equation}\label{bound_for_E_v_NEW}
{\rm \textbf{Pr}}(E_v)<\frac{2}{\Delta^3}\leq \frac{1}{\Delta^2\cdot\Delta^\frac{2}{3}}
\end{equation}
for $\Delta$ sufficiently large.

Now consider \emph{(VI)}.
For a vertex $v\in V$ of degree $d$, let finally
$F_v$ be the event that  $d\geq \Delta/3$ and for some integer $\alpha\in\{1,2,\ldots,\lceil 3\Delta^{\frac{1}{3}}/\ln^\frac{1}{3}\Delta\rceil\}$,
the number of neighbours $u$ of $v$ with $d(u)\geq \Delta/3$ and
$S(u)\in I_{\alpha}$ is greater than
$\Delta^{\frac{5}{6}}\ln^{\frac{1}{6}}\Delta+ \Delta^{\frac{1}{2}}$
(note that for $\alpha \geq \lceil 3\Delta^{\frac{1}{3}}/\ln^\frac{1}{3}\Delta\rceil+1$, $I_{\alpha}\subset (\Delta^2,+\infty)$,
while $S(u)$ are assumed not exceed $\Delta^2$).
For a given vertex $v$ of degree $d\geq \Delta/3$ and an integer $\alpha\in\{1,2,\ldots,\lceil 3\Delta^{\frac{1}{3}}/\ln^\frac{1}{3}\Delta\rceil\}$,
let then $F_{v,\alpha}$
denote the random variable of the number of neighbours
$u\in N(v)$
with $d(u)\geq \Delta/3$ and $S(u)\in I_{\alpha}$.
For every fixed neighbour $u$ of $v$ with $d(u)\geq \Delta/3$, by~(\ref{I_d_alpha}) and (\ref{S_v_in_lemma}) we have:
\begin{eqnarray}
{\rm \textbf{Pr}}\left(S(u) \in I_{\alpha}\right)&=&
{\rm \textbf{Pr}}\Bigg(\left(\left\lceil\Delta^{\frac{2}{3}}\ln^{\frac{1}{3}}\Delta\right\rceil
+6\left\lceil\Delta^{\frac{1}{3}}\ln^{\frac{2}{3}}\Delta\right\rceil\right)d(u)c_1(u)\label{S_u_equality}\\
&&\in \left((\alpha'-1) \frac{\Delta^{\frac{5}{3}}\ln^{\frac{1}{3}}\Delta}{3}, \alpha'\frac{\Delta^{\frac{5}{3}}\ln^{\frac{1}{3}}\Delta}{3}\right]\Bigg)\nonumber
\end{eqnarray}
for some real number $\alpha'$. Since
$\frac{\frac{\Delta^{\frac{5}{3}}\ln^{\frac{1}{3}}\Delta}{3}}
{(\lceil\Delta^{\frac{2}{3}}\ln^{\frac{1}{3}}\Delta\rceil
+6\lceil\Delta^{\frac{1}{3}}\ln^{\frac{2}{3}}\Delta\rceil)d(u)}
\leq 1$,
the probability
in~(\ref{S_u_equality}) may by bounded from above by the following one,
where $\alpha''$ is also some constant real number,
$${\rm \textbf{Pr}}\left(S(u) \in I_{\alpha}\right) \leq {\rm \textbf{Pr}}\left(c_1(u) \in (\alpha''-1,\alpha'']\right)
\leq \frac{1}{\lceil\Delta^{\frac{1}{6}}\ln^{-\frac{1}{6}}\Delta\rceil}\leq\frac{1}{\Delta^{\frac{1}{6}}\ln^{-\frac{1}{6}}\Delta}.$$
Note that since the value of $S(u)$ depends only on the choice of $c_1(u)$,
then in our random process, $S(u)$, for $u\in V$, are \emph{independent} random variables.
Therefore, by the Chernoff Bound,
\begin{eqnarray*}
{\rm \textbf{Pr}}\left(F_{v,\alpha} >\Delta^{\frac{5}{6}}\ln^{\frac{1}{6}}\Delta+ \Delta^{\frac{1}{2}}\right)
&\leq& {\rm \textbf{Pr}}\left(BIN\left(\Delta,\frac{1}{\Delta^{\frac{1}{6}}\ln^{-\frac{1}{6}}\Delta}\right)
>\Delta^{\frac{5}{6}}\ln^{\frac{1}{6}}\Delta+ \Delta^{\frac{1}{2}}\right)\\
&\leq& {\rm \textbf{Pr}}\left(\left|BIN\left(\Delta,\frac{1}{\Delta^{\frac{1}{6}}\ln^{-\frac{1}{6}}\Delta}\right)
-\Delta^{\frac{5}{6}}\ln^{\frac{1}{6}}\Delta\right|>\Delta^{\frac{1}{2}}\right)\\
&<& 2e^{-\frac{\Delta}{3\Delta^{\frac{5}{6}}\ln^{\frac{1}{6}}\Delta}}
\leq 2e^{-3\ln\Delta} = \frac{2}{\Delta^3},
\end{eqnarray*}
and hence (for $d\geq \Delta/3$),
\begin{eqnarray}
{\rm \textbf{Pr}}(F_v)&\leq& \sum_{\alpha=1}^{\left\lceil\frac{3\Delta^{\frac{1}{3}}}{\ln^\frac{1}{3}\Delta}\right\rceil}
{\rm \textbf{Pr}}\left(F_{v,\alpha} >\Delta^{\frac{5}{6}}\ln^{\frac{1}{6}}\Delta+ \Delta^{\frac{1}{2}}\right)
< \left\lceil\frac{3\Delta^{\frac{1}{3}}}{\ln^\frac{1}{3}\Delta}\right\rceil\cdot \frac{2}{\Delta^3}
\leq \frac{1}{\Delta^2\cdot\Delta^\frac{2}{3}}.\label{bound_for_F_v_NEW}
\end{eqnarray}

Note that since each of the events $A_v$, $B_v$, $C_v$, $D_v$, $E_v$ and $F_v$ depends only on the
random colours of $v$ and its adjacent vertices or edges,
then each such event corresponding to a vertex $v$ is mutually independent of all other events corresponding to vertices $v'$ at distance at least three from $v$,
hence is mutually independent of all except at most $D=5+6\Delta^2$ other events. Moreover, by (\ref{bound_for_A_v}),  (\ref{bound_for_B_v}),  (\ref{bound_for_C_v_NEW}), (\ref{bound_for_D_v_NEW}), (\ref{bound_for_E_v_NEW}) and  (\ref{bound_for_F_v_NEW}),
the probability of each of these events equals at most
$\frac{1}{\Delta^2\cdot\Delta^\frac{2}{3}}$.
Since
$$e\frac{1}{\Delta^2\cdot\Delta^\frac{2}{3}}(6+6\Delta^2) < 1,$$
then by the Lov\'asz Local Lemma we thus obtain:
$${\rm \textbf{Pr}}\left(\bigcap_{v\in V}\overline{A_v}\cap \overline{B_v}\cap \overline{C_v}\cap \overline{D_v}\cap \overline{E_v}\cap \overline{F_v}\right)>0,$$
and therefore there exist colourings
$c_1:V\to\{1,2,\ldots,\lceil\Delta^\frac{1}{6}/\ln^\frac{1}{6}\Delta\rceil\}$,
$c_2:E\to\{1,2,\ldots,$ $\lceil\Delta^{\frac{1}{3}}/\ln^\frac{1}{3}\Delta\rceil\}$ and $c_3:V\to\{1,2,\ldots,\lceil\Delta^\frac{1}{3}/\ln^\frac{1}{3}\Delta\rceil\}$
satisfying $(I),(II),(VI)$, $(1^\circ)$ and $(2^\circ)$.

In the following second step of our construction, we shall first define (a preliminary) colouring
$c_3:E\to \{1,2,\ldots,2\lceil\Delta^{\frac{1}{6}}/\ln^\frac{1}{6}\Delta\rceil +\lceil\Delta^{\frac{1}{3}}/\ln^\frac{1}{3}\Delta\rceil\}$
(or more formally -- extend the vertex colouring $c_3$ into the total colouring
$c_3:V\cup E\to \{1,2,\ldots,2\lceil\Delta^{\frac{1}{6}}/\ln^\frac{1}{6}\Delta\rceil +\lceil\Delta^{\frac{1}{3}}/\ln^\frac{1}{3}\Delta\rceil\}$)
and then exploit once more the Lov\'asz Local Lemma coupled with the Chernoff's theorem on concentration
to make small rearrangements
of the $c_3(e)$'s for a limited number of edges $e\in E$,
aimed at satisfying \emph{(III)-(V)}.
No changes shall be introduced to $c_1(v)$, $c_2(e)$, $c_3(v)$ with $v\in V$, $e\in E$ though,
and thus our requirements \emph{(I), (II)} and \emph{(VI)} shall remain fulfilled.

Initially, let us set $c_3(uv)=c_1(u)+c_1(v)+c_2(uv)$ for \emph{all} edges $uv\in E$.
Let $H_1=(V,E_1)$ be a subgraph of $G$, where $uv\in E_1$ if and only if
$c_1(u)+c_1(v)+c_2(uv)=c_3(u)$ or $c_1(u)+c_1(v)+c_2(uv)=c_3(v)$,
and let $H_2=(V,E_2)$ be a subgraph of $G$, where $uv\in E_2$ if and only if
$c_3(u)=c_3(v)$.
Note that by $(1^\circ)$, $\Delta(H_1)\leq 2\Delta^\frac{2}{3}\ln^\frac{1}{3}\Delta+5\Delta^\frac{1}{3}\ln^\frac{2}{3}\Delta$,
and by $(2^\circ)$, $\Delta(H_2)\leq \Delta^\frac{2}{3}\ln^\frac{1}{3}\Delta+3\Delta^\frac{1}{3}\ln^\frac{2}{3}\Delta$.
Further on we shall only modify the colours $c_3(e)$ for the edges $e\in E_1\cup E_2$,
thus the requirement \emph{(III)} shall be fulfilled at the end of the construction.

Let us uncolour all edges $e\in E_1$, e.g., by temporarily setting $c_3(e)=0$ for these.
Then let us set (recolour) $c_3(uv)=c_3(u)$ ($=c_3(v)$) for every edge $uv\in E_2$.
Note that at this point, by $(1^\circ)$ and $(2^\circ)$, we have that:
\begin{itemize}
\item[($3^\circ$)]for every $v\in V$
and every $c_3^*\in\{1,2,\ldots,2\lceil\Delta^{\frac{1}{6}}/\ln^\frac{1}{6}\Delta\rceil+\lceil\Delta^{\frac{1}{3}}/\ln^\frac{1}{3}\Delta\rceil\}$,
the number of edges
$uv$ (incident with $v$) for which $c_3(uv)=c_3^*$
does not exceed
$\Delta^\frac{2}{3}\ln^\frac{1}{3}\Delta + 3\Delta^\frac{1}{3}\ln^\frac{2}{3}\Delta$.
\end{itemize}
We then choose a new colouring for the remaining uncoloured edges.
Denote the subgraph they induce by $H_3$, thus $H_3=(V,E_1\smallsetminus E_2)$ and
hence $\Delta(H_3)\leq 2\Delta^\frac{2}{3}\ln^\frac{1}{3}\Delta+5\Delta^\frac{1}{3}\ln^\frac{2}{3}\Delta$.
For every $e\in E_1\smallsetminus E_2$ assign randomly, independently and equiprobably a new colour
$c_3(e)\in
\{1,2,\ldots,2\lceil\Delta^{\frac{1}{6}}/\ln^\frac{1}{6}\Delta\rceil+\lceil\Delta^{\frac{1}{3}}/\ln^\frac{1}{3}\Delta\rceil\}$.     We claim that this can be done so that afterwards:
\begin{itemize}
\item[($4^\circ$)] for every $v\in V$
and every $c_3^*\in\{1,2,\ldots,2\lceil\Delta^{\frac{1}{6}}/\ln^\frac{1}{6}\Delta\rceil+\lceil\Delta^{\frac{1}{3}}/\ln^\frac{1}{3}\Delta\rceil\}$,
the number of edges
$uv\in E_1\smallsetminus E_2$ for which $c_3(uv)=c_3^*$
does not exceed $2\Delta^\frac{1}{3}\ln^\frac{2}{3}\Delta$.
\end{itemize}
To see that, for any given $v\in V$, let $L_v$ be the event that
the number of edges
$uv\in E_1\smallsetminus E_2$ with
$c_3(uv)=c_3^*$ exceeds $2\Delta^\frac{1}{3}\ln^\frac{2}{3}\Delta$
for some $c_3^*\in\{1,\ldots,2\lceil\Delta^{\frac{1}{6}}/\ln^\frac{1}{6}\Delta\rceil+\lceil\Delta^{\frac{1}{3}}/\ln^\frac{1}{3}\Delta\rceil\}$.
For every $c_3^*\in\{1,\ldots,2\lceil\Delta^{\frac{1}{6}}/\ln^\frac{1}{6}\Delta\rceil+\lceil\Delta^{\frac{1}{3}}/\ln^\frac{1}{3}\Delta\rceil\}$,
 denote further by $J_{v,c_3^*}$ the random variable of the number of edges
$uv\in E_1\smallsetminus E_2$ with $c_3(uv)=c_3^*$.
Note that
\begin{equation}\label{auxiliary_ineq_JP}
\frac{1}
{\lceil\Delta^\frac{1}{3}\ln^{-\frac{1}{3}}\Delta\rceil+2\lceil\Delta^\frac{1}{6}\ln^{-\frac{1}{6}}\Delta\rceil}
\leq
\frac{2\Delta^\frac{1}{3}\ln^\frac{2}{3}\Delta-3\Delta^\frac{1}{6}\ln^\frac{5}{6}\Delta}
{\lfloor 2\Delta^\frac{2}{3}\ln^\frac{1}{3}\Delta+5\Delta^\frac{1}{3}\ln^\frac{2}{3}\Delta\rfloor}.
\end{equation}
By the Chernoff Bound and (\ref{auxiliary_ineq_JP}) we thus obtain (for sufficiently large $\Delta$):
\begin{eqnarray*}
&&{\rm \textbf{Pr}}\left(J_{v,c_3^*} > 2\Delta^\frac{1}{3}\ln^\frac{2}{3}\Delta \right)\\
&\leq &
{\rm \textbf{Pr}}
\Bigg(BIN\left(\lfloor2\Delta^\frac{2}{3}\ln^\frac{1}{3}\Delta+5\Delta^\frac{1}{3}\ln^\frac{2}{3}\Delta\rfloor,
\frac{1}{\lceil\Delta^\frac{1}{3}\ln^{-\frac{1}{3}}\Delta\rceil+2\lceil\Delta^\frac{1}{6}\ln^{-\frac{1}{6}}\Delta\rceil}\right)
>2\Delta^\frac{1}{3}\ln^\frac{2}{3}\Delta\Bigg)\\
&\leq &
{\rm \textbf{Pr}}
\Bigg(BIN\left(\lfloor2\Delta^\frac{2}{3}\ln^\frac{1}{3}\Delta+5\Delta^\frac{1}{3}\ln^\frac{2}{3}\Delta\rfloor,
\frac{2\Delta^\frac{1}{3}\ln^\frac{2}{3}\Delta-3\Delta^\frac{1}{6}\ln^\frac{5}{6}\Delta}
{\lfloor 2\Delta^\frac{2}{3}\ln^\frac{1}{3}\Delta+5\Delta^\frac{1}{3}\ln^\frac{2}{3}\Delta\rfloor}\right)\\
&&>\left(2\Delta^{\frac{1}{3}}\ln^{\frac{2}{3}}\Delta- 3\Delta^{\frac{1}{6}}\ln^{\frac{5}{6}}\Delta\right)
+3\Delta^{\frac{1}{6}}\ln^{\frac{5}{6}}\Delta\Bigg)\\
&\leq &
{\rm \textbf{Pr}}
\Bigg(\Bigg|BIN\left(\lfloor2\Delta^\frac{2}{3}\ln^\frac{1}{3}\Delta+5\Delta^\frac{1}{3}\ln^\frac{2}{3}\Delta\rfloor,
\frac{2\Delta^\frac{1}{3}\ln^\frac{2}{3}\Delta-3\Delta^\frac{1}{6}\ln^\frac{5}{6}\Delta}
{\lfloor 2\Delta^\frac{2}{3}\ln^\frac{1}{3}\Delta+5\Delta^\frac{1}{3}\ln^\frac{2}{3}\Delta\rfloor}\right)\\
&&-\left(2\Delta^{\frac{1}{3}}\ln^{\frac{2}{3}}\Delta- 3\Delta^{\frac{1}{6}}\ln^{\frac{5}{6}}\Delta\right)\Bigg|
>3\Delta^{\frac{1}{6}}\ln^{\frac{5}{6}}\Delta\Bigg)\\
&<& 2e^{-\frac{9\Delta^{\frac{1}{3}}\ln^{\frac{5}{3}}\Delta}
{3(2\Delta^{\frac{1}{3}}\ln^{\frac{2}{3}}\Delta- 3\Delta^{\frac{1}{6}}\ln^{\frac{5}{6}}\Delta)}}
< 2e^{-\frac{9\Delta^{\frac{1}{3}}\ln^{\frac{5}{3}}\Delta}
{6\Delta^{\frac{1}{3}}\ln^{\frac{2}{3}}\Delta}}
= \frac{2}{\Delta^\frac{3}{2}}.
\end{eqnarray*}
Consequently,
\begin{equation}\label{bound_for_L_v_NEW}
{\rm \textbf{Pr}}(L_v) \leq \left(2\left\lceil\frac{\Delta^{\frac{1}{6}}}{\ln^\frac{1}{6}\Delta}\right\rceil
+\left\lceil\frac{\Delta^{\frac{1}{3}}}{\ln^\frac{1}{3}\Delta}\right\rceil\right) \cdot  \frac{2}{\Delta^\frac{3}{2}} < \frac{1}{\Delta}.
\end{equation}
Note that since every event $L_v$ depends only on the random colours
assigned to the edges incident with $v$ in $H_3$,
then each such event is mutually independent of all
other events $L_{v'}$ with $v'$ at distance greater than one from $v$,
i.e., of all except at most $\Delta(H_3)$ other events.
Since
$$e\frac{1}{\Delta}\left(\Delta(H_3)+1\right)\leq
e\frac{1}{\Delta}\left(2\Delta^\frac{2}{3}\ln^\frac{1}{3}\Delta+5\Delta^\frac{1}{3}\ln^\frac{2}{3}\Delta+1\right)
<1,$$
by (\ref{bound_for_L_v_NEW}) and the Lov\'asz Local Lemma we thus obtain that
$${\rm \textbf{Pr}}\left(\bigcap_{v\in V}\overline{L_v}\right)>0,$$
hence our claim ($4^\circ$) holds.
By ($3^\circ$) and ($4^\circ$), the obtained colouring $c_3$ is consistent with the requirement \emph{(IV)}
from the thesis.
Note   also that in our construction, the colours $c_3(e)$
for the edges of the subgraph $H_2$ were finally defined accordingly to the requirement \emph{(V)}.
\qed
\end{pf}

\section{Main Result}
\begin{theorem}\label{main_JP_Th_as}
If $G$ is a graph
of sufficiently large maximum degree $\Delta$, then
\begin{equation}\label{main_JP_ineq}
\chi''_{\sum}(G) \leq \Delta+139\Delta^{\frac{5}{6}}\ln^\frac{1}{6}\Delta,
\end{equation}
i.e., $\chi''_{\sum}(G)\leq(1+o(1))\Delta$.
\end{theorem}
Within our approach, we
made effort to optimizing the order of the second term of the upper bound (\ref{main_JP_ineq}).
With more
scrupulous calculus,
one may decrease the constant $139$ though.

\section{Proof of the Main Result}\label{section_with_main_result}
\subsection{Preliminary Edge Colouring}
Let $G=(V,E)$ be a graph.
Whenever needed we shall
assume that its maximum degree $\Delta$ is sufficiently large.

Let $c_1:V\to \{1,2,\ldots,\lceil\Delta^\frac{1}{6}/\ln^\frac{1}{6}\Delta\rceil\}$,
$c_2: E\to \{1,2,\ldots,$ $\lceil\Delta^\frac{1}{3}/\ln^\frac{1}{3}\Delta\rceil\}$ and
$c_3:V\cup E \to\{1,2,\ldots,2\lceil\Delta^{\frac{1}{6}}/\ln^\frac{1}{6}\Delta\rceil +\lceil\Delta^{\frac{1}{3}}/\ln^\frac{1}{3}\Delta\rceil\}$
be auxiliary vertex, edge and total, resp,. colourings of $G$
guaranteed by Lemma~\ref{main_probabilistic_lemma}, where the function $R(d,\Delta)$ shall be specified later (and may be derived from~(\ref{S_v}) below).

Set $B= (\lceil\Delta^{\frac{2}{3}}\ln^\frac{1}{3}\Delta\rceil
+6\lceil\Delta^{\frac{1}{3}}\ln^\frac{2}{3} \Delta\rceil)$.
For every edge $e\in E$ and every vertex $v\in V$,
set their initial \emph{temporary colours} as
\begin{equation}\label{definition_of_c'}
c_t(e)=B\cdot c_3(e), ~~~~c_t(v)=B\cdot c_3(v),
\end{equation}
hence,
\begin{eqnarray}\label{c_prime_bounds}
c_t(e),c_t(v)\leq \Delta+2\Delta^{\frac{5}{6}}\ln^\frac{1}{6}\Delta + o(\Delta^{\frac{5}{6}}\ln^\frac{1}{6}\Delta)
\end{eqnarray}
for every $e\in E$, $v\in V$.
These colours shall be modified in the further part of the construction,
while $c_t(e),c_t(v)$ shall always refer to the up-to-date edge and vertex colours,
non of which
shall ever exceed
$\Delta+139\Delta^{\frac{5}{6}}\ln^\frac{1}{6}\Delta$.

\subsection{Making the Total Colouring Proper}
Let $C$ be the constant from Corollary~\ref{MolloReedCor},
and assume that $\Delta$ is large enough so that
$C\leq \Delta^{\frac{1}{3}}\ln^\frac{2}{3}\Delta$.

For every $\beta\in \{1,2,\ldots,$
$2\lceil\Delta^{\frac{1}{6}}/\ln^\frac{1}{6}\Delta\rceil$ $+\lceil\Delta^{\frac{1}{3}}/$ $\ln^\frac{1}{3}\Delta\rceil\}$, the set of integers
$$
\{(\beta-1) B+1,(\beta-1) B+2,\ldots,\beta B\}
$$
shall be called a \emph{colour class} (corresponding to $\beta$).
We shall also abuse this notation and write that an edge
or a vertex with a colour from this set \emph{belongs} to this colour class.
Note that at this point of the construction, only the largest
integers from the colour classes
might appear as
colours of the edges and vertices of $G$.
We shall first modify these preliminary colours
in order to make the total colouring of
$G$ proper,
but at the same time, these colours shall remain in unchanged colour classes
throughout (almost) the whole argument, until the beginning of its final paragraph.
Thus if $c_t(e)=Bc_3(e)+a(e)$ or $c_t(v)=Bc_3(v)+a(v)$
in any following, except the last stage of the construction, then
\begin{equation}\label{range_of_a_e}
a(e),a(v)\in\left\{-B+1,-B+2,\ldots,0\right\}.
\end{equation}
The quantity $a(e)$ (or $a(v)$), which we shall call the \emph{addition},
simply determines which element from the colour class associated with $e$ (or $v$)
is assigned to this edge (or vertex) at a given stage of the process.
Note that by~(\ref{c_prime_bounds}) this also guarantees that
we will further on have (almost until the end of the proof):
\begin{equation}\label{c_prime_upper_bound}
c_t(e),c_t(v) \leq \Delta+(2+o(1))\Delta^{\frac{5}{6}}\ln^\frac{1}{6}\Delta
\end{equation}
for every $e\in E$, $v\in V$.

By Lemma~\ref{main_probabilistic_lemma}\emph{(IV)} (and (\ref{definition_of_c'})),
for every colour class and a vertex $v$,
the number of its incident edges belonging to this colour class
does not exceed
$\Delta^{\frac{2}{3}}\ln^{\frac{1}{3}}\Delta+5\Delta^{\frac{1}{3}}\ln^\frac{2}{3}\Delta$,
hence the maximum degree of every subgraphs induced by the edges from any single
colour class in $G$ does not exceed this quantity as well.
Since $C\leq \Delta^{\frac{1}{3}}\ln^\frac{2}{3}\Delta$,
by Corollary~\ref{MolloReedCor}, one by one
for every such subgraph we may find the additions
$a(e),a(v)\in\{-B+1,-B+2,\ldots,0\}$
for all edges and vertices so that the obtained (temporary) total colouring of $G$ is proper,
where every vertex retains only colour chosen for it while fixing the total colouring
for the subgraph induced by the edges from its colour class.
This way, by Lemma~\ref{main_probabilistic_lemma}\emph{(V)} and our construction,
for every edge $uv\in E$, $c_t(u)\neq c_t(v)$,
hence the total colouring obtained
is indeed proper.

\subsection{Adjustment of Weighted Degrees}
Now, in order to
make our colouring neighbour sum distinguishing,
we have to first show that the weighted degrees of the neighbours of every vertex $v$
(of large degree)
are indeed well distributed. Let us first note that at this point of the construction,
if $v$ is a vertex of degree
$d\geq\Delta/3$,
then its weighted degree is the following:
$$
s(v)=c_t(v)+\sum_{u\in N(v)}c_t(uv)
=\left(Bc_3(v)+a(v)\right)+\sum_{u\in N(v)}\left(Bc_3(uv)+a(uv)\right).
$$
By Lemma~\ref{main_probabilistic_lemma}\emph{(III)}, we thus may write that:
\begin{equation}
s(v)
= \sum_{u\in N(v)}\left(B\left[c_1(u)+c_1(v)+c_2(uv)\right]
+a(uv)\right)+f(v),\label{first_s_v_equality}
\end{equation}
with (since edge colours are all in the range $1,\ldots,\Delta+O(\Delta^\frac{5}{6}\ln^\frac{1}{6}\Delta)$):
\begin{equation}\label{modul_fv_ineq}
|f(v)|\leq (3+o(1))\Delta^\frac{2}{3}\ln^\frac{1}{3}\Delta \left(\Delta+O(\Delta^\frac{5}{6}\ln^\frac{1}{6}\Delta)\right)
\leq (3+o(1))\Delta^\frac{5}{3}\ln^\frac{1}{3}\Delta
\end{equation}
(where the quantity $Bc_3(v)+a(v)$ is counted in $f(v)$).
Moreover, by~(\ref{range_of_a_e}):
\begin{equation}\label{a_e_sum}
-\Delta B \leq \sum_{u\in N(v)}a(uv) \leq 0,
\end{equation}
by Lemma~\ref{main_probabilistic_lemma}\emph{(I)}:
\begin{eqnarray}\sum_{u\in N(v)}c_1(u) &=& \sum_{i=1}^{\lceil\Delta^{\frac{1}{6}}\ln^{-\frac{1}{6}}\Delta\rceil}
\left(\frac{d}{\lceil\Delta^{\frac{1}{6}}\ln^{-\frac{1}{6}}\Delta\rceil}+ f_{1,i}(v)\right) i\nonumber\\
& =& \left(\frac{d}{\lceil\Delta^{\frac{1}{6}}\ln^{-\frac{1}{6}}\Delta\rceil} + f_1(v)\right)
{\lceil\Delta^{\frac{1}{6}}\ln^{-\frac{1}{6}}\Delta\rceil+1 \choose 2},\label{c_1_sum}
\end{eqnarray}
where $f_1(v)$ is some real number with
\begin{equation}\label{f_1_bound}
|f_1(v)|\leq \Delta^\frac{1}{2}
\end{equation}
(see Lemma~\ref{main_probabilistic_lemma}\emph{(I)}), and by Lemma~\ref{main_probabilistic_lemma}\emph{(II)}:
\begin{eqnarray}\sum_{u\in N(v)}c_2(uv) &=& \sum_{i=1}^{\lceil\Delta^{\frac{1}{3}}\ln^{-\frac{1}{3}}\Delta\rceil}
\left(\frac{d}{\lceil\Delta^{\frac{1}{3}}\ln^{-\frac{1}{3}}\Delta\rceil}+ f_{2,i}(v)\right)i\nonumber\\
& =& \left(\frac{d}{\lceil\Delta^{\frac{1}{3}}\ln^{-\frac{1}{3}}\Delta\rceil} + f_2(v) \right)
{\lceil\Delta^{\frac{1}{3}}\ln^{-\frac{1}{3}}\Delta\rceil+1 \choose 2},\label{c_2_sum}
\end{eqnarray}
where $f_2(v)$ is some real number with
\begin{equation}\label{f_2_bound}
|f_2(v)|\leq 3\Delta^{\frac{1}{3}}\ln^{\frac{2}{3}}\Delta
\end{equation}
(see Lemma~\ref{main_probabilistic_lemma}\emph{(II)}).

By~(\ref{first_s_v_equality}), (\ref{c_1_sum}) and (\ref{c_2_sum}) we thus obtain that:
$$s(v)=S(v)+ F(v),$$
where
\begin{equation}\label{S_v}
S(v)=B
\left[dc_1(v)+\frac{d}{\lceil\Delta^{\frac{1}{6}}\ln^{-\frac{1}{6}}\Delta\rceil}
{\lceil\Delta^{\frac{1}{6}}\ln^{-\frac{1}{6}}\Delta\rceil+1 \choose 2}
+\frac{d}{\lceil\Delta^{\frac{1}{3}}\ln^{-\frac{1}{3}}\Delta\rceil}{\lceil\Delta^{\frac{1}{3}}\ln^{-\frac{1}{3}}\Delta\rceil+1 \choose 2}\right]
\end{equation}
is the function required to apply Lemma~\ref{main_probabilistic_lemma} (with $S(v)\leq \Delta^2$ for $\Delta$ sufficiently large), and
\begin{eqnarray*}
F(v)&=&f(v)+\sum_{u\in N(v)}a(uv) +
B\left(f_1(v)
{\lceil\Delta^{\frac{1}{6}}\ln^{-\frac{1}{6}}\Delta\rceil+1 \choose 2} +
f_2(v)
{\lceil\Delta^{\frac{1}{3}}\ln^{-\frac{1}{3}}\Delta\rceil+1 \choose 2}\right),
\end{eqnarray*}
for which, by~(\ref{modul_fv_ineq}), (\ref{a_e_sum}), (\ref{f_1_bound}) and (\ref{f_2_bound}), the following holds:
\begin{eqnarray}
|F(v)|&\leq& (3+o(1))\Delta^\frac{5}{3}\ln^\frac{1}{3}\Delta\nonumber\\
&&+ \left(\left\lceil\Delta^{\frac{2}{3}}\ln^\frac{1}{3}\Delta\right\rceil
+6\left\lceil\Delta^{\frac{1}{3}}\ln^\frac{2}{3} \Delta\right\rceil\right)\nonumber\\
&&\times
\Bigg(\Delta + {\lceil\Delta^{\frac{1}{6}}\ln^{-\frac{1}{6}}\Delta\rceil+1 \choose 2}
\Delta^\frac{1}{2}
+3{\lceil\Delta^{\frac{1}{3}}\ln^{-\frac{1}{3}}\Delta\rceil+1 \choose 2} \Delta^{\frac{1}{3}}\ln^{\frac{2}{3}}\Delta\Bigg)\nonumber\\
&=& \left(\frac{11}{2}+o(1)\right)\Delta^\frac{5}{3}\ln^\frac{1}{3}\Delta.\label{F_v_second_bound}
\end{eqnarray}
This \emph{fault} $F(v)$ will not change significantly
(i.e., it shall remain true that $|F(v)|\leq (\frac{11}{2}+o(1))\Delta^\frac{5}{3}\ln^\frac{1}{3}\Delta$)
even if we later arbitrarily change the colours of no more than $15\ln\Delta$
edges incident with $v$ (each by at most $\Delta+O(\Delta^\frac{5}{6}\ln^\frac{1}{6}\Delta)$),
hence the weighted degree of $v$ by less than $16\Delta\ln\Delta$.

By our construction, for every vertex $v$ of large degree,
only some part of its neighbours (of large degree)
might at this point have their temporary weighted degrees
close enough to the one of $v$ to be threatened with a potential conflict with $v$.
We shall estimate an upper bound for the number of such neighbours.

Consider a vertex $v$ of degree $d\geq \Delta/3$
and its neighbour $u$
with $d(u)\geq \Delta/3$.
Then by (\ref{F_v_second_bound}),
\begin{eqnarray*}
&&\left(F(v)+16\Delta\ln\Delta \right)
+\left(F(u)+16\Delta\ln\Delta \right)\\
&\leq& \left(11+o(1)\right)\Delta^\frac{5}{3}\ln^\frac{1}{3}\Delta
= \left(33+o(1)\right)\frac{\Delta^{\frac{5}{3}} \ln^\frac{1}{3}\Delta}{3},
\end{eqnarray*}
where $\Delta^{\frac{5}{3}} \ln^\frac{1}{3}\Delta/3$ is the length of
the intervals $I_{\alpha}$ (cf. (\ref{I_d_alpha})).
Therefore, $s(u)$ might eventually be equal to $s(v)$ only if $S(u)$ belongs
to a certain interval, say $K_v$, of length
$2(33+o(1))\Delta^{\frac{5}{3}} \ln^\frac{1}{3}\Delta/3 = (66+o(1))\Delta^{\frac{5}{3}} \ln^\frac{1}{3}\Delta/3$ centered in $S(v)$.
Since such $K_v$ may be incident with at most $68$ intervals $I_\alpha$,
if we denote by $R_v$ the set of all neighbours of $v$ with this (risky) feature
(i.e., with $d(u)\geq \Delta/3$ and $S(u)\in K_v$),
then by Lemma~\ref{main_probabilistic_lemma}\emph{(VI)} we thus have:
\begin{equation}
|R_v|\leq 68 \left(\Delta^{\frac{5}{6}}\ln^{\frac{1}{6}}\Delta
+\Delta^{\frac{1}{2}}\right)
= \left(68+o(1)\right)\Delta^\frac{5}{6}\ln^\frac{1}{6}\Delta.\label{R_v_bound}
\end{equation}

Finally we shall use a relatively short list of new colours
to recolour few of the edges so that all neighbours of degree $\geq\Delta/3$
have pairwise distinct weighted degrees in $G$. Distinguishing the remaining vertices
shall be straightforward than.
First we choose a subgraph of $G$
to which every vertex of large degree (in $G$) contributes
at least two of its incident edges.

For every vertex $v$ of degree $d\geq \Delta/3$,
we independently choose a pair of its incident edges randomly and equiprobably, and we denote
the graph induced by the set of the chosen edges (some of these possibly twice) by $H$.
Note that
$d_H(v)\geq 2$ for such vertices. On the other hand,
given a vertex $v\in V$ of degree $d$, the probability that an
edge $uv\in E$ was chosen to $H$ `by a neighbour' $u$ of $v$
equals at most $6/\Delta$ (since this may only have happened if $d(u)\geq\Delta/3$,
and moreover, each of the incident edges of $u$ might belong to the pair chosen for $u$
with probability $2/d(u)$), where the choices of all neighbours are independent.
Taking into account two more possible edges incident with $v$ in $H$ (those chosen by $v$ in the case when $d\geq \Delta/3$),
by the Chernoff Bound:
\begin{eqnarray*}
{\rm \textbf{Pr}}\left(d_H(v)-2>14\ln\Delta\right) &\leq&
{\rm \textbf{Pr}}\left(BIN\left(d,\frac{6}{\Delta}\right)>14\ln\Delta\right)\\
&\leq& {\rm \textbf{Pr}}\left(BIN\left(\Delta,\frac{7\ln\Delta}{\Delta}\right)>14\ln\Delta\right)\\
&\leq& {\rm \textbf{Pr}}\left(\left|BIN\left(\Delta,\frac{7\ln\Delta}{\Delta}\right)-7\ln\Delta\right|
>7\ln\Delta\right)\\
&<&2e^{-\frac{7\ln\Delta}{3}} = \frac{2}{\Delta^\frac{7}{3}}.
\end{eqnarray*}
Note that if we denote by $W_v$ the event that $d_H(v)-2>14\ln\Delta$, then
every $W_v$ is mutually independent of all events $W_{v'}$ with $v'$ at distance at least three from $v$,
i.e., of all but at most $D=\Delta^2$ (other) events.
 By the Lov\'asz Local Lemma we thus may choose our subgraph $H$ so that
\begin{equation}\label{d_H_v-1-bound}
d_{H}(v)\leq 14\ln\Delta +2\leq 15\ln\Delta
\end{equation}
for every $v\in V$.

Recall that thus far for every edge $e\in E$ and vertex $v\in V$,
$c_t(e),c_t(v)\leq \Delta+(2+o(1))\Delta^{\frac{5}{6}}\ln^\frac{1}{6}\Delta$,
see~(\ref{c_prime_upper_bound}).
We shall now use integer colours not exceeding
$\Delta+139\Delta^\frac{5}{6}\ln^\frac{1}{6}\Delta$
which were not yet used for any edge nor any vertex of $G$.
Denote the set of these by $A$ and note that by~(\ref{c_prime_upper_bound}),
(\ref{R_v_bound}) and (\ref{d_H_v-1-bound}),
for every $u,v\in V$ with $d(u),d(v)\geq\Delta/3$,
\begin{equation}\label{A_vs_R_u_R_v}
|A|\geq|R_u|+|R_v|+2\Delta(H)
\end{equation}
for sufficiently large $\Delta$.
One by one, we analyze subsequent edges of $H$ in any order,
and for every consecutive edge $uv\in E(H)$
we recolour $uv$ using a colour in $A$ not yet used by any of the (at most $2(\Delta(H)-1)$)
edges incident with $uv$ in $H$
so that $u$ obtains the weighted degree distinct from all its neighbours in $R_u$ (if $d(u)\geq \Delta/3$),
except possibly $v$,
and $v$ obtains the weighted degree distinct from all its neighbours in $R_v$ (if $d(v)\geq \Delta/3$),
except possibly $u$.
By~(\ref{A_vs_R_u_R_v}) this is always feasible.
The obtained
colours $c(e)=c_t(e)$ of the edges $e\in E$
make up their final colouring.
Note that since $d_H(v)\geq 2$ for every vertex $v$ of degree $\geq \Delta/3$,
at this stage, every vertex of degree at least $\Delta/3$
has a weighted degree distinct from \emph{all} its neighbours of degree at least $\Delta/3$.
To eliminate the remaining conflicts, we adjust the colours of vertices of smaller degrees.
Thus one after another, for every vertex $v\in V$ of degree $d(v)<\Delta/3$ we simply
choose a (possibly new) colour from $\{1,2,\ldots,\Delta\}$ for $v$ so that the weighted degree of $v$
is distinct from those of its neighbours
and the total colouring remains proper.
Since $v$ is adjacent with less than $\Delta/3$ coloured vertices (and edges)
with weighted degrees to be avoided, there is always an available choice of colour for $v$.
The obtained vertex colours $c(v)$ are the final ones,
and the resulting proper total colouring $c$ of $G$ is neighbour sum distinguishing.
\qed

\section{Final remarks}
In fact
one may obtain the same result as in Theorem~\ref{main_JP_Th_as}
using a slightly shorter, but very similar argument.
For this aim, it is sufficient to apply Kahn's result~\cite{Kahn}
 on the renowned List Colouring Conjecture,
rather than Theorem~\ref{MolloyReedTh} (or Corollary~\ref{MolloReedCor}),
preceded by a few simplifications in the main probabilistic lemma
and in the colouring procedure itself.
We have however decided to apply the (likely more sharp) result of Molloy and Reed
on the Total Colouring Conjecture, as it seems
more promising in view of possible
further improvements of our result.
In particular it would be interesting to prove the following
weaker version of Conjecture~\ref{PilsniakWozniakCon}.
\begin{conjecture}
There exists a constant $C$ such that $\chi''_{\sum}(G)\leq \Delta(G)+C$ for every graph $G$.
\end{conjecture}

Curiously, neither of these two approaches is useful while trying
to obtain a result concerning the natural list version of the problem investigated
(where colours are being chosen from arbitrary lists of real number of possibly shortest lengths),
not even the one based on the application of Kahn's result on list colourings.
It thus would also be very interesting to invent a new method
resolving the following problem.
\begin{conjecture}
The natural list correspondent of $\chi''_{\sum}(G)$ is bounded from above by
$(1+o(1))\Delta$.
\end{conjecture}

\end{document}